\documentclass[11pt]{amsart}%
\usepackage{amssymb,amsmath,amsfonts,latexsym,amsthm,geometry,graphicx}
\usepackage{amsmath}
\usepackage{amsfonts}
\usepackage{amssymb}
\usepackage{graphicx}%
\providecommand{\U}[1]{\protect\rule{.1in}{.1in}}
\providecommand{\U}[1]{\protect\rule{.1in}{.1in}}
\providecommand{\U}[1]{\protect\rule{.1in}{.1in}}
\providecommand{\U}[1]{\protect\rule{.1in}{.1in}}
\newtheorem{theorem}{Theorem}[section]

\theoremstyle{definition}

\newtheorem{remark}[theorem]{Remark}

\begin{document}
\title[Bohnenblust--Hille inequality]{Bohnenblust--Hille inequality for polynomials whose monomials have uniformly
bounded number of variables}
\author[M. Maia]{Mariana Maia}
\address[M. Maia]{Departamento de Matem\'{a}tica\\
\indent Universidade Federal da Para\'{\i}ba\\
\indent 58.051-900 - Jo\~{a}o Pessoa, Brazil.}
\email{mariana.britomaia@gmail.com}
\author[T. Nogueira]{Tony Nogueira}
\address[T. Nogueira]{Departamento de Matem\'{a}tica\\
\indent Universidade Federal da Para\'{\i}ba\\
\indent 58.051-900 - Jo\~{a}o Pessoa, Brazil.}
\email{tonykleverson@gmail.com}
\author[D. Pellegrino]{Daniel Pellegrino}
\address[D. Pellegrino]{Departamento de Matem\'{a}tica \\
Universidade Federal da Para\'{\i}ba \\
58.051-900 - Jo\~{a}o Pessoa, Brazil.}
\email{pellegrino@pq.cnpq.br and dmpellegrino@gmail.com}
\thanks{M. Maia and T. Nogueira are supported by Capes and D. Pellegrino is supported
by CNPq.}
\keywords{Bohnenblust--Hille inequality}
\subjclass[2010]{}

\begin{abstract}
In 2015, using an innovative technique, Carando, Defant and Sevilla-Peris
succeeded in proving a Bohnenblust--Hille type inequality with constants of
polynomial growth in $m$ for a certain family of complex $m$-homogeneous
polynomials. In the present paper, using a completely different approach, we
prove that the constants of this inequality are uniformly bounded
irrespectively of the value of $m$.

\end{abstract}
\maketitle

\section{Introduction}

The Bohnenblust--Hille inequality \cite{bh} for complex $m$-homogeneous
polynomials asserts that there is a constant $C_{m}>0$ such that
\[
\left(  \sum_{\left\vert \alpha\right\vert =m}|c_{\alpha}(P)|^{\frac{2m}{m+1}%
}\right)  ^{\frac{m+1}{2m}}\leq C_{m}\Vert P\Vert
\]
for all continuous $m$-homogeneous polynomials $P:c_{0}\rightarrow\mathbb{C}$
of the form $P(x)=\sum_{|\alpha|=m}c_{\alpha}(P)\mathbf{x}^{\alpha}$. This
inequality is important in many fields of Mathematics. In 2011, it has been
proven in \cite{ann} that $C_{m}$ can be chosen with exponential growth, and
this result had several important applications. In 2014 the estimates of
\cite{ann} were improved in \cite{bohr} and it has been shown that for any
$\varepsilon>0$ there is a constant $\kappa>0$ such that%
\[
C_{m}\leq\kappa\left(  1+\varepsilon\right)  ^{m}%
\]
and this result was crucial to obtain the final solution to the asymptotic
growth of the Bohr radius problem. In 2015, Carando, Defant and Sevilla-Peris
have shown that for a particular family of $m$-homogeneous polynomials the
constants $C_{m}$ could have been chosen to have polynomial growth in
$m.$\ More precisely, they have proved that for polynomials whose monomials
have a uniformly bounded number $M$ of different variables, there is a
Bohnenblust-Hille type inequality with a constant of polynomial growth in $m$.
Our main result shows, by means of a completely different technique, that in
fact these constants are uniformly bounded by a constant that does not depend
on $m.$

Let $\mathbb{K}=\mathbb{R}$ or $\mathbb{C}$, let $\alpha=(\alpha_{j}%
)_{j=1}^{\infty}$ be a sequence in $\mathbb{N}\cup\{0\}$ and$,$ as usual,
define $|\alpha|=%
{\textstyle\sum}
\alpha_{j};$ in this case we also denote $\mathbf{x}^{\alpha}:=%
{\textstyle\prod\nolimits_{j}}
x_{j}^{\alpha_{j}}.$ An $m$-homogeneous polynomial $P:c_{0}\rightarrow
\mathbb{C}$ is denoted by
\[
P(x)=\sum_{|\alpha|=m}c_{\alpha}(P)\mathbf{x}^{\alpha}.
\]
We recall that the norm of $P$ is given by $\Vert P\Vert:=\sup_{x\in B_{c_{0}%
}}|P(x)|$. Since $|\alpha|=m$, it is obvious that only a finite number of
$\alpha_{j}$ are non null and we define $\binom{m}{\alpha}:=\frac{m!}%
{\alpha_{1}!\dots\alpha_{n}!},$ where $\alpha_{1},...,\alpha_{n}$ are the non
null elements of $\alpha$. Following the notation of \cite{carando}, for
positive integers $m$ and $M\leq m$ we define
\[
vars(P)=\text{card }\{j:\alpha_{j}\neq0\}
\]
and
\[
\Lambda_{M}=\{\alpha:|\alpha|=m,\quad vars(P)\leq M\}.
\]
In \cite{carando} it was proved that%
\begin{equation}
\left(  \sum_{\alpha\in\Lambda_{M}}|c_{\alpha}(P)|^{\frac{2m}{m+1}}\right)
^{\frac{m+1}{2m}}\leq2^{\frac{M}{2}}m^{\frac{M+1}{2}}\Vert P\Vert\label{ju}%
\end{equation}
for all continuous $m$-homogeneous polynomials $P:c_{0}\rightarrow\mathbb{C}$,
being stressed by the authors of \cite{carando} the polynomial growth (in $m$)
of the constants -- this was a very nice property having in mind that in
general the best known estimates had just a subexponential growth (in $m$).
Our main result shows that the optimal constants of (\ref{ju}) are universally
bounded irrespectively of the value of $m$:

\begin{theorem}
\label{arrocho} For all positive integers $m$ and $M\leq m$, there exists a
constant $\kappa_{M}>0$ such that
\[
\left(  \sum_{\alpha\in\Lambda_{M}}|c_{\alpha}(P)|^{\frac{2m}{m+1}}\right)
^{\frac{m+1}{2m}}\leq\kappa_{M}\Vert P\Vert
\]
for all continuous $m$-homogeneous polynomials $P:c_{0}\rightarrow\mathbb{C}$.
\end{theorem}

\section{The proof of Theorem \ref{arrocho}}

The main result of \ \cite{www} asserts that if $1\leq k\leq m$ and
$n_{1},\dots,n_{k}\geq1$ are positive integers such that $n_{1}+\cdots
+n_{k}=m$, then there is a constant $C_{k,m}^{\mathbb{K}}\geq1$ such that
\begin{equation}
\left(  \sum_{i_{1},\dots,i_{k}=1}^{\infty}\left\vert T\left(  e_{i_{1}%
}^{n_{1}},\dots,e_{i_{k}}^{n_{k}}\right)  \right\vert ^{\frac{2k}{k+1}%
}\right)  ^{\frac{k+1}{2k}}\leq C_{k,m}^{\mathbb{K}}\left\Vert T\right\Vert
\label{clas}%
\end{equation}
for all continuous $m$--linear forms $T:c_{0}\times\dots\times c_{0}%
\rightarrow\mathbb{K}$. Moreover, the exponent $\frac{2k}{k+1}$ is optimal. In
\cite{www} it is also proved that
\begin{equation}
C_{k,m}^{\mathbb{K}}\leq C_{k}^{\mathbb{K}} \label{c22}%
\end{equation}
for all $1\leq k\leq m$, where $C_{k}^{\mathbb{K}}$ is the optimal constant of
the $k$-linear Bohnenblust--Hille inequality. From \cite{bohr} we know that
there are positive constants $\beta_{1},\beta_{2}$ such that
\begin{equation}
C_{k}^{\mathbb{C}}\leq\beta_{1}k^{\frac{1-\gamma}{2}}<\beta_{1}k^{0.212}
\label{9ll}%
\end{equation}
and%
\[
C_{k}^{\mathbb{R}}\leq\beta_{2}k^{\frac{2-\log2-\gamma}{2}}<\beta_{2}%
k^{0.365},
\]
where $\gamma$ is the Euler-Mascheroni constant. From now on $C_{k}%
^{\mathbb{K}}$ will be just denoted by $C_{k}.$

Let $\hat{P}$ be the symmetric $m$-linear form associated to $P$. Let%
\[
\Gamma_{m}:=\left\{  \tau=\left(  \tau_{1},...,\tau_{M}\right)  \in
\lbrack0,m]^{M}:\tau_{1}+\dots+\tau_{M}=m\right\}  .
\]
It is simple to verify that (an argument of symmetry provides even sharper
estimates but, surprisingly, this rough estimate is enough for our purposes)
\[
\sum_{\alpha\in\Lambda_{M}}|c_{\alpha}(P)|^{\frac{2M}{M+1}}\leq\sum_{\tau
\in\Gamma_{m}}\sum_{i_{1},\dots,i_{M}}\binom{m}{\tau}^{\frac{2M}{M+1}}|\hat
{P}(e_{i_{1}}^{\tau_{1}},\dots,e_{i_{M}}^{\tau_{M}})|^{\frac{2M}{M+1}}.
\]
If $\lfloor x\rfloor:=\max\{n\in\mathbb{N}:n\leq x\}$ it is a simple exercise
to verify that
\begin{equation}
\binom{m}{\tau}\leq\frac{m!}{\left(  \lfloor\frac{m}{M}\rfloor!\right)  ^{M}%
}.\label{des1}%
\end{equation}
By \eqref{des1}, since
\[
\sum_{i_{1},\dots,i_{k}}|\hat{P}(e_{i_{1}}^{\tau_{1}},\dots,e_{i_{k}}%
^{\tau_{k}})|^{\frac{2M}{M+1}}\leq\sum_{i_{1},\dots,i_{M}}|\hat{P}(e_{i_{1}%
}^{\tau_{1}},\dots,e_{i_{M}}^{\tau_{M}})|^{\frac{2M}{M+1}}%
\]
for all $1\leq k\leq M,$ we have
\begin{align*}
\sum_{i_{1},\dots,i_{M}}\binom{m}{\tau}^{\frac{2M}{M+1}}|\hat{P}(e_{i_{1}%
}^{\tau_{1}},\dots,e_{i_{m}}^{\tau_{M}})|^{\frac{2M}{M+1}} &  \leq\left(
\frac{m!}{\left(  \lfloor\frac{m}{M}\rfloor!\right)  ^{M}}\right)  ^{\frac
{2M}{M+1}}\sum_{i_{1},\dots,i_{M}}|\hat{P}(e_{i_{1}}^{\tau_{1}},\dots
,e_{i_{m}}^{\tau_{M}})|^{\frac{2M}{M+1}}\\
&  \leq\left(  \frac{m!}{\left(  \lfloor\frac{m}{M}\rfloor!\right)  ^{M}%
}\right)  ^{\frac{2M}{M+1}}(C_{M}\Vert\hat{P}\Vert)^{\frac{2M}{M+1}},
\end{align*}
where in the last inequality we have used (\ref{clas}) and (\ref{c22}). Note
also that
\begin{align*}
\sum_{\alpha\in\Lambda_{M}}|c_{\alpha}(P)|^{\frac{2M}{M+1}} &  \leq\sum
_{\tau\in\Gamma_{M}}\left(  \frac{m!}{\left(  \lfloor\frac{m}{M}%
\rfloor!\right)  ^{M}}\right)  ^{\frac{2M}{M+1}}(C_{M}\Vert\hat{P}%
\Vert)^{\frac{2M}{M+1}}\\
&  =\binom{m+M-1}{m}\left(  \frac{m!}{\left(  \lfloor\frac{m}{M}%
\rfloor!\right)  ^{M}}\right)  ^{\frac{2M}{M+1}}(C_{M}\Vert\hat{P}%
\Vert)^{\frac{2M}{M+1}}.
\end{align*}
Hence
\begin{align*}
\left(  \sum_{\alpha\in\Lambda_{M}}|c_{\alpha}(P)|^{\frac{2M}{M+1}}\right)
^{\frac{M+1}{2M}} &  \leq\binom{m+M-1}{m}^{\frac{M+1}{2M}}\frac{m!}{\left(
\lfloor\frac{m}{M}\rfloor!\right)  ^{M}}C_{M}\Vert\hat{P}\Vert\\
&  \leq\binom{m+M-1}{m}^{\frac{M+1}{2M}}\frac{m!}{\left(  \lfloor\frac{m}%
{M}\rfloor!\right)  ^{M}}C_{M}e^{m}\Vert P\Vert.
\end{align*}
The above estimate holds for complex and real scalars. Since now we are just
dealing with complex scalars, by the Maximum Modulus Principle, we have
\[
\left(  \sum_{\alpha\in\Lambda_{M}}|c_{\alpha}(P)|^{2}\right)  ^{\frac{1}{2}%
}\leq\left(  \sum_{|\alpha|=m}|c_{\alpha}(P)|^{2}\right)  ^{\frac{1}{2}}%
\leq\Vert P\Vert.
\]
Since%
\[
\frac{1}{\frac{2m}{m+1}}=\frac{\theta}{\frac{2M}{M+1}}+\frac{1-\theta}{2}%
\]
with
\[
\theta=\frac{M}{m},
\]
by a corollary of the H\"{o}lder inequality, we have
\begin{align}
\left(  \sum_{\alpha\in\Lambda_{M}}|c_{\alpha}(P)|^{\frac{2m}{m+1}}\right)
^{\frac{m+1}{2m}} &  \leq\left[  \left(  \sum_{\alpha\in\Lambda_{M}}%
|a_{\alpha}|^{\frac{2M}{M+1}}\right)  ^{\frac{M+1}{2M}}\right]  ^{\frac{M}{m}%
}\left[  \left(  \sum_{\alpha\in\Lambda_{M}}|a_{\alpha}|^{2}\right)
^{\frac{1}{2}}\right]  ^{1-\frac{M}{m}}\label{999}\\
&  \leq\left(  \binom{m+M-1}{m}^{\frac{M+1}{2M}}\frac{m!}{\left(  \lfloor
\frac{m}{M}\rfloor!\right)  ^{M}}C_{M}e^{m}\Vert P\Vert\right)  ^{\frac{M}{m}%
}\Vert P\Vert^{1-\frac{M}{m}}\nonumber
\end{align}
Using the Stirling Formula we can prove that%
\begin{equation}
\lim_{m\rightarrow\infty}\left[  \frac{m!}{\left(  \lfloor\frac{m}{M}%
\rfloor!\right)  ^{M}}\right]  ^{\frac{M}{m}}=M^{M}.\label{0uy}%
\end{equation}
By (\ref{0uy}) and (\ref{999}) we conclude that there is a constant $\zeta
_{M}>0$ such that%
\begin{align*}
\left(  \sum_{\alpha\in\Lambda_{M}}|c_{\alpha}(P)|^{\frac{2m}{m+1}}\right)
^{\frac{m+1}{2m}} &  \leq\binom{m+M-1}{m}^{\frac{M+1}{2m}}\left[  \frac
{m!}{\left(  \lfloor\frac{m}{M}\rfloor!\right)  ^{M}}\right]  ^{\frac{M}{m}%
}C_{M}^{\frac{M}{m}}e^{M}\Vert P\Vert\\
&  \leq\zeta_{M}\binom{m+M-1}{m}^{\frac{M+1}{2m}}C_{M}^{\frac{M}{m}}e^{M}\Vert
P\Vert
\end{align*}
Recalling that $\gamma$ denotes the Euler--Mascheroni constant, by the
previous inequality combined with (\ref{9ll}) we have%
\begin{align*}
\left(  \sum_{\alpha\in\Lambda_{M}}|c_{\alpha}(P)|^{\frac{2m}{m+1}}\right)
^{\frac{m+1}{2m}} &  \leq\zeta_{M}\binom{m+M-1}{m}^{\frac{M+1}{2m}}\left(
\beta_{1}M^{\frac{1-\gamma}{2}}\right)  ^{\frac{M}{m}}e^{M}\Vert P\Vert\\
&  \leq\kappa_{M}\Vert P\Vert
\end{align*}
for a certain $\kappa_{M}>0$ and this concludes the proof.

\begin{remark}
By the above proof we also conclude that for this Bohnenblust--Hille type
inequality the optimal exponent is not $\frac{2m}{m+1}$, contrary to what
happens for the classical Bohnenblust--Hille inequality. We can see that the
inequality holds for the smaller exponent $\frac{2M}{M+1}$ and using the
optimality of the exponent $\frac{2M}{M+1}$ in \cite{www} it seems to be just
an exercise to show that this exponent $\frac{2M}{M+1}$ is sharp in this case.


\end{remark}

\end{document}